\documentclass[aoas,preprint]{imsart}

\RequirePackage[OT1]{fontenc}
\RequirePackage{amsthm,amsmath}
\RequirePackage[numbers]{natbib}
\RequirePackage[colorlinks,citecolor=blue,urlcolor=blue]{hyperref}
\usepackage[pdftex]{graphicx}

\startlocaldefs
\numberwithin{equation}{section}
\theoremstyle{plain}
\newtheorem{Theorem}{Theorem}[section]

\endlocaldefs

\begin{document}

\begin{frontmatter}
\title{Multi-Armed Bandit Problem\\
and Batch UCB Rule} \runtitle{Multi-Armed Bandit Problem and Batch
UCB Rule}

\begin{aug}
\author{\fnms{Alexander} \snm{Kolnogorov}\thanksref{m1}\ead[label=e1]{Alexander.Kolnogorov@novsu.ru}}
\and
\author{\fnms{Sergey} \snm{Garbar}\thanksref{m1}
\ead[label=e2]{Sergey.Garbar@novsu.ru} }

\runauthor{A. Kolnogorov et al.}

\affiliation{Yaroslav-the-Wise Novgorod State
University\thanksmark{m1}}

\address{41 B.Saint-Petersburgskaya Str.,
Velikiy Novgorod, Russia, 173003\\
Applied Mathematics and Information Science Department\\
\printead{e1}\\
\printead{e2}}

\end{aug}

\begin{abstract}
We obtain the upper bound of the loss function for a strategy in
the multi-armed bandit problem with Gaussian distributions of
incomes. Considered strategy is an asymptotic generalization of
the strategy proposed by J. Bather for the multi-armed bandit
problem and using UCB rule, i.e. choosing the action corresponding
to the maximum of the upper bound of the confidence interval of
the current estimate of the expected value of one-step income.
Results are obtained with the help of invariant description of the
control on the unit horizon in the domain of close distributions
because just there the loss function attains its maximal values.
UCB rule is widely used in machine learning. It can be also used
for the batch data processing optimization if there are two
alternative processing methods available with different a priori
unknown efficiencies.
\end{abstract}



\end{frontmatter}

\section{Introduction}

We consider the multi-armed bandit problem.  Multi-armed bandit is
a slot machine with two or more arms the choice of which is
accompanied with a random income of a gambler  depending only on
chosen arm ~\cite{Berry}. The goal of a gambler is to maximize the
total expected income. To this end, meanwhile the game process, he
should determine the arm corresponding to the largest expected
one-step income and provide its predominant usage. The problem is
also well-known as the problem of expedient
behavior~\cite{Tsetlin} and of adaptive control in a random
environment~\cite{Sragovich}. The problem has numerous
applications in machine learning~\cite{Auer, Lugosi}.

In what follows, we consider Gaussian multi-armed bandit which
naturally arises when batch data processing is optimized and there
are two or more processing methods available with different a
priori unknown efficiencies~\cite{Koln}. Formally, it is a
controlled random process $\xi_n$, $n=1,2,\dots,N$, which value at
the point of time $n$ depends only on currently chosen arm $y_n$,
is interpreted as income and has Gaussian (normal) distribution
with probability density $f_D(x|m_\ell)=(2\pi
D)^{-1/2}\exp\left(-(x-m_\ell)^2/(2D)\right)$ if $y_n=\ell$,
$\ell=1,\dots,J$. Variance $D$ is assumed to be known and
expectations $m_1,\dots,m_J$ are assumed to be unknown. The
requirement of a priori known variance can be omitted because
considered algorithm only a little changes under a moderate change
of the  variance (e.g., 5--10\% change), hence, the variance can
be estimated at the initial stage of the control.

A control strategy $\sigma$ determines, generally, a randomized
choice of action $y_n$ depending on currently available
information of the history of the process. In what follows, we
restrict consideration with the following strategy proposed
in~\cite{Bather}. Let the $\ell$-th action be applied $n_\ell$
times up to the point of time $n$ and let $X_\ell$ denote
corresponding cumulative income ($\ell=1,\dots,J$). In this case
$X_\ell/n_\ell$ is a current estimate of the expectation $m_\ell$.
Since the goal is to maximize the total expected income, it seems
reasonable always to apply the action corresponding to currently
largest value $X_\ell/n_\ell$. However, it is well known that such
a rule can result in a significant losses because by chance the
initial estimate $X_\ell/n_\ell$, corresponding to the largest
$m_\ell$, can take a low value and, therefore, this action will be
never applied in the sequel. Instead of estimates of $\{m_\ell\}$
themselves let's consider the upper bounds of their confidence
intervals
\begin{gather}\label{kolnogorovAV1}
    U_\ell(n)=\frac{X_\ell(n)}{n_\ell}+\frac{aD^{1/2}}{n_\ell^{1/2}}(2+\zeta_\ell(n)),
\end{gather}
where $a>0$, $\{\zeta_l(n)\}$ are i.i.d. random variables with
probability density $e^{-x}$ ($x>0$); $\ell=1,2,\dots,J$;
$n=1,2,\dots,N$.

Considered strategy prescribes initially once to apply all the
actions by turns and then at each point of time $n$ to choose the
action corresponding to the largest value $\{U_\ell(n)\}$.
Strategies of such a form are called UCB (Upper Confidence Bound)
rules. Considered in this paper strategy at $a=2/15$ is equivalent
to the strategy proposed in~\cite{Bather} for Bernoulli
multi-armed bandit problem (up to summands of the order
$n_\ell^{-1}$). In this case one should put $D=0.25$ which is
equal to the maximum value of the variance of Bernoulli one-step
income. It is noted in~\cite{Bather} that at $J=2$ the maximal
expected losses (scaled to the value $(DN)^{1/2}$) do not exceed
0.72 for large $N$. However, explanation of this result is not
presented in~\cite{Bather} and, to the best of our knowledge, it
was not published later.

We explain this result using the invariant description of the
control on the unit horizon in the domain of close distributions
where the maximum values of expected losses are attained in the
multi-armed bandit problem. Moreover, we consider the batch
version of the strategy~\cite{Bather} and show that expected
losses depend only on the number of processed batches and some
invariant characteristics of the parameter. Note that batch
(parallel) strategies are especially important when processing
time of the data item is significant, because in this case the
total processing time depends on the number of batches rather than
on the total number of data. The maximum scaled expected losses
for the batch UCB rule are turned out to be 0.75, i.e. are almost
the same as in~\cite{Bather}. We also note that different versions
of the UCB rule are widely used in machine learning (see,
e.g.~\cite{Auer, Lugosi}).

\section{Main Results}

Considered multi-armed bandit can be described by a vector
parameter  $\theta=(m_1,\dots,m_J)$. Let's define the loss
function. If the parameter is known, one should always choose the
action corresponding to the maximum of $m_1,\dots,m_J$, the total
expected income would thus be equal to $N\max(m_1,\dots,m_J)$. For
actually applied strategy $\sigma$ the total expected income is
less than maximal possible by the value which is called the loss
function and is equal to
\begin{gather}\label{kolnogorovAV2}
    L_N(\sigma,\theta)=\mathbf{E}_{\sigma,\theta}\left(\sum_{n=1}^N \left(\max(m_1,\dots,m_J)-
   \xi_n \right)\right).
\end{gather}
Here $\mathbf{E}_{\sigma,\theta}$ denotes mathematical expectation
calculated with respect to measure generated by strategy $\sigma$
and parameter $\theta$. We are interested in the upper bound of
the maximum losses calculated on the set of admissible values of
parameter which is chosen the following
\begin{gather*}
\Theta=\{m_\ell=m+d_\ell (D/N)^{1/2}; m \in (-\infty,+\infty),
|d_\ell| \le C < \infty, \ell=1,\dots,J\}.
\end{gather*}
This is the set of parameters describing ``close'' distributions
which are characterized by the difference of mathematical
expectations of the order $N^{-1/2}$. Maximal expected losses are
attained just there and have the order  $N^{1/2}$ (see, e.g,
\cite{Vogel}). For ``distant'' distributions the losses have
smaller values. For example, they are of the order $\ln(N)$ if
$\max(m_1,\dots,m_J)$ exceeds all other $\{m_\ell\}$ by some
$\delta>0$ (see, e.g., \cite{Lai}).

Let's consider strategies which change the actions only after
applying them $M$ times in succession. These strategies allow
batch (and also parallel) processing. We assume for simplicity
that $N=MK$ where $K$ is the number of batches. For batch
strategies the upper bounds~\eqref{kolnogorovAV1} take the form
\begin{gather}\label{kolnogorovAV3}
    U_\ell(k)=\frac{X_\ell(k)}{k_\ell}+\frac{a(MD)^{1/2}}{k_\ell^{1/2}}(2+\zeta_\ell(k)),
\end{gather}
where $k$ is the number of batches and $k_\ell$, $X_\ell(k)$ are
the cumulative number of batches to which the $\ell$-th action was
applied and corresponding cumulative income after processing  $k$
batches ($k=1,2,\dots,K$). Let's denote by
\begin{gather*}
    I_\ell(k)=\left\{
    \begin{array}{l}
    1, \quad \mbox{if } U_\ell(k)=\max(U_1(k),\dots,U_J(k)),\\
    0, \quad \mbox{otherwise}
     \end{array}\right.
\end{gather*}
the indicator of chosen action for processing the $(k+1)$-th batch
according to considered rule at $k>J$ (recall that at $k\le J$
actions are chosen by turns). Note that with probability 1 only
one value of $\{I_\ell(k)\}$ is equal to 1. For considered
parameter the following presentation holds
\begin{gather}\label{kolnogorovAV4}
X_\ell(k)= k_\ell M \left(m+d_\ell
\left(\frac{D}{N}\right)^{1/2}\right)+ \sum_{i=1}^k I_\ell(i)
Y_\ell (MD;i),
\end{gather}
where $\{Y_\ell (MD;i)\}$ are i.i.d. normally distributed random
variables with zero mathematical expectations and variances equal
to $MD$. Let's introduce the variables $t=kK^{-1}$, $t_\ell=k_\ell
K^{-1}$, $\varepsilon=K^{-1}$. By \eqref{kolnogorovAV3},
\eqref{kolnogorovAV4}, it follows that
\begin{gather*}
    U_\ell(k)=   Mm+d_l\left(\frac{MD}{K}\right)^{1/2}+\\+\frac{(MD)^{1/2}
   \sum_{i=1}^k  I_\ell(i)
    Y_\ell (\varepsilon;i)}{t_\ell K^{1/2}}+\frac{a(MD)^{1/2}}{t_\ell^{1/2} K^{1/2}}(2+\zeta_\ell(n)),
\end{gather*}
$\ell=1,2,\dots,J$; $k=J+1,J+2,\dots,K$. After the linear
transformation $u_\ell(t)=(U_\ell(k)-Mm)(MD)^{-1/2}K^{1/2}$, which
does not change the arrangement of bounds, we obtain the upper
bounds in invariant form with a control horizon equal to
 1:
\begin{gather}\label{kolnogorovAV5}
    u_\ell(t)=d_l+\frac{\sum_{i=1}^k I_\ell(i)  Y_\ell (\varepsilon;i)}{t_\ell}+\frac{a}{t_\ell^{1/2} }(2+\zeta_\ell(t)),
\end{gather}
$\ell=1,2,\dots,J$; $t=(J+1)\varepsilon,(J+2)\varepsilon,\dots,1$.

Let's determine the loss function. For chosen parameter without
loss of generality let's assume that $d_1=\max(d_1,\dots,d_J)$.
Then
\begin{gather*}
L_N(\sigma,\theta)=(D/N)^{1/2}\sum_{\ell=2}^J
(d_1-d_\ell)\mathbf{E}_{\sigma,\theta}\left(\sum_{k=1}^{K}M
I_\ell(k)\right)=\\= (DN)^{1/2}\sum_{\ell=2}^K
(d_1-d_\ell)\mathbf{E}_{\sigma,\theta}\left(\sum_{k=1}^{K}\varepsilon
I_\ell(k)\right),
\end{gather*}
and for scaled (by $(DN)^{1/2}$) loss function we draw the
expression
\begin{gather}\label{kolnogorovAV6}
(DN)^{-1/2} L_N(\sigma,\theta)=\sum_{\ell=2}^K
(d_1-d_\ell)\mathbf{E}_{\sigma,\theta}\left(\sum_{k=1}^{K}\varepsilon
I_\ell(k)\right).
\end{gather}

\begin{figure}[htb]
\begin{center}
\includegraphics[scale=0.65]{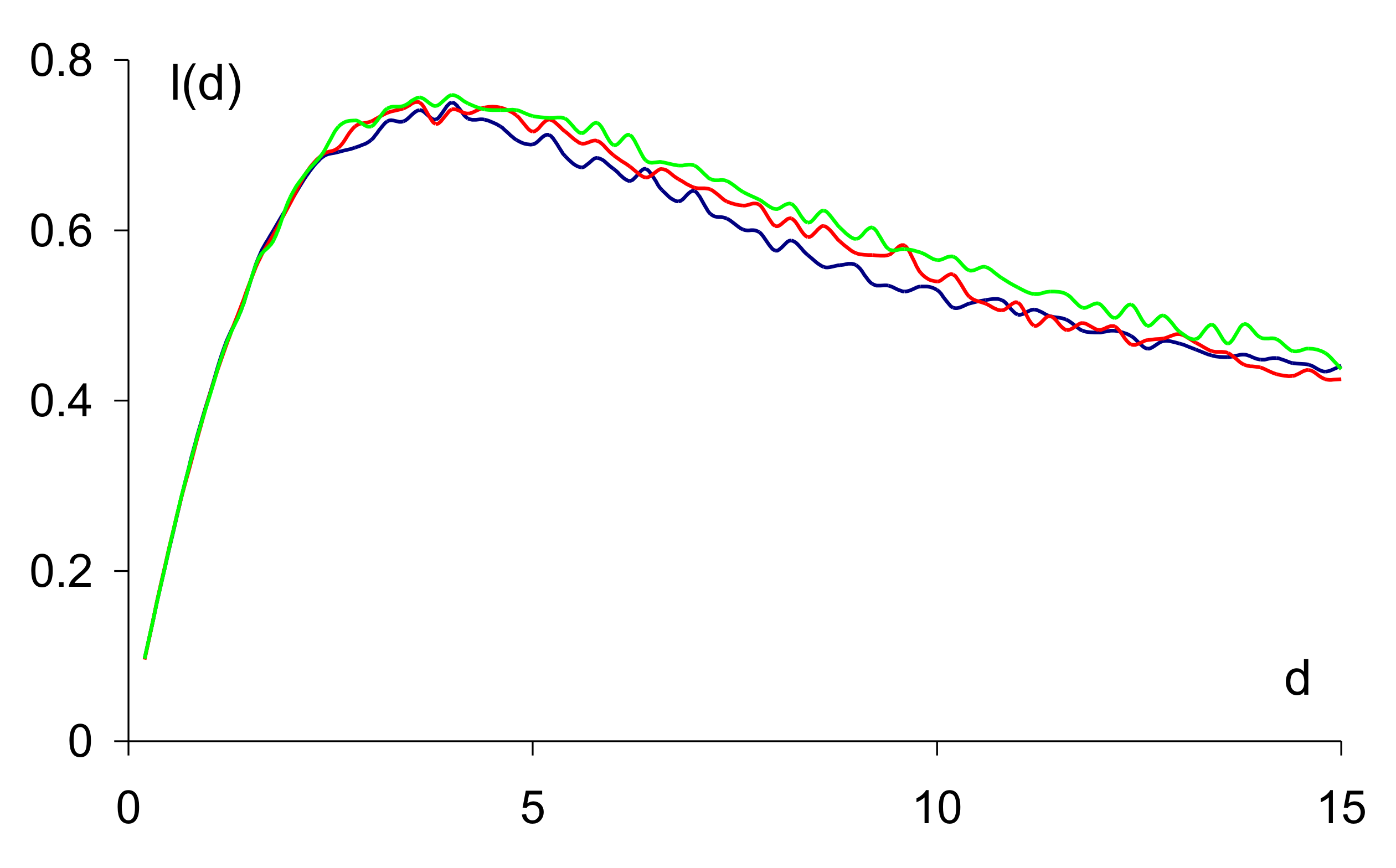}
\caption{Scaled Expected Losses}\label{fig1}
\end{center}
\end{figure}

Results can be presented by the following theorem.

\begin{Theorem}
The usage of the rule~\eqref{kolnogorovAV3} results in invariant
description on the unit control horizon which is described
by~\eqref{kolnogorovAV5}. For scaled (by $(DN)^{1/2}$) loss
function the expression~\eqref{kolnogorovAV6} holds. Expressions
\eqref{kolnogorovAV5}, \eqref{kolnogorovAV6} depend only on the
number of processed batches.
\end{Theorem}

On figure Fig.~\ref{fig1} we present Monte-Carlo simulation
results for the scaled loss function corresponding to UCB rule at
$a=1/3$ and $K=2$. Averaging was implemented over 10000
simulations. Given $K=2$, one can take $d_1=-d_2=0.5d$,
corresponding scaled loss function is denoted by $l(d)$. Blue, red
and green lines are obtained for control horizons
$N=100,400,1500$. One can see that $\max l(d)\approx 0.75$ at $d
\approx 3.5$. Hence, the maximum value of the scaled loss function
is close to determined in~\cite{Bather}. But the value $a=1/3$
which is close to optimal in considered case, i.e. providing the
minimum of maximal losses, significantly differs from 2/15
determined in~\cite{Bather}. This may be explained by the fact
that in~\cite{Bather} considerably small control horizons were
considered, e.g., $N=50$. In this case the summands of the order
$n_\ell^{-1}$, which are present in the rule proposed
in~\cite{Bather}, essentially affect the values $\{U_\ell(n)\}$.
%%%%%%%%%%%%%%%%%%%%%%%%%%%%%%%%%%%%%%%%%%
%\vspace{6pt}

\end{document}